\documentclass [12pt,a4paper,reqno]{amsart}
\usepackage{amsmath,amsthm}
\usepackage{amsfonts}
\usepackage{amssymb}
\allowdisplaybreaks

\newcommand{\be}{\begin{equation}}
\newcommand{\ef}{\end{equation}}
\chardef\bslash=`\\ 

\newtheorem*{thm*}{Theorem}

\theoremstyle{definition}

\newtheorem*{remark*}{Remarks}
\newtheorem*{defn*}{Definition}

\theoremstyle{remark}

\newcommand{\G}{\Gamma}
\newcommand{\wt}{\widetilde}
\newcommand{\wh}{\widehat}
\newcommand{\fc}{\frac}

\newcommand{\iy}{\infty}
 \renewcommand{\sectionmark}[1]{}

\newcommand{\Be}{Beltrami}

\newcommand{\qc} {quasiconformal}

\newcommand{\ve}{\varepsilon}

\newcommand{\Te} {Teichm\"{u}ller}

 \usepackage{amsfonts}
\newcommand{\field}[1]{\mathbb{#1}}
\newcommand{\g}{\gamma}

\newcommand{\D}{\field{D}}

\newcommand{\ov}{\overline}
\newcommand{\vp}{\varphi}
\newcommand{\hC}{\wh{\field{C}}}
\newcommand{\C}{\field{C}}

\newcommand{\B}{\mathbf{B}}
\newcommand{\T}{\mathbf{T}}
\newcommand{\Belt}{\operatorname{Belt}}

\newcommand{\Mob}{\operatorname{Mob}}

\newcommand{\Teich}{\operatorname{Teich}}

\newcommand{\x} {\mathbf x}
\renewcommand{\a} {\alpha}

\newcommand{\ld}{\lambda}
\newcommand{\kp}{\kappa}

\begin{document}

\title{Teichm\"{u}ller space theory and classical problems of geometric function theory}
\author{Samuel L. Krushkal}

\begin{abstract} Recently the author presented a new approach to solving the coefficient problems for holomorphic functions based on the deep features of Teichm\"{u}ller spaces.It involves the Bers isomorphism theorem for Teichm\"{u}ller spaces of punctured Riemann surfaces.

The aim of the present paper is to provide new applications of this approach and extend the indicated results to more general classes of functions

\end{abstract}

\date{\today\hskip4mm({TeichSpCoefProbl.tex})}

\maketitle

\bigskip

{\small {\textbf {2010 Mathematics Subject Classification:} Primary: 30C50, 30C75, 30F60;
Secondary 30C55, 30C62, 31A05, 32L05, 32Q45}

\medskip

}

\bigskip

\markboth{Samuel L. Krushkal}{Teichm\"{u}ller space theory and classical
problems} \pagestyle{headings}

\bigskip
{\small\em{ \centerline{To Vladimir Gutlyanskii on the occasion of his 80th birthday}}}

\bigskip\bigskip
\centerline{\bf 1. INTRODUCTION, CLASSES OF FUNCTIONS AND STATEMENT}
\centerline{\bf OF MAIN RESULTS}

\bigskip\noindent
Recently the author presented in \cite{Kr5} a new approach
to solving the classical coefficient problems on various classes of holomorphic
functions, not necessarily univalent. This approach involves the deep analytic and
geometric features of Teichm\"{u}ller spaces, especially the Bers isomorphism
theorem for Teichm\"{u}ller spaces of punctured Riemann surfaces.

Estimating holomorphic functionals depending on the Taylor coefficients
of univalent holomorphic functions is important in various geometric and physical applications of complex analysis.

In the present paper, we provide new applications of this approach and extend the results of \cite{Kr5} to more general classes of functions.

\bigskip\noindent
{\bf 1.1}. We start with the general collection $\wh S(1)$ of univalent functions on the unit disk $\D = \{|z| < 1\}$ which
is the completion in the topology of locally uniform convergence on $\D$ of the
set of univalent functions
$$
w(z) = a_1 z + a_2 z^2 + \dots \quad \text{with} \ \ |a_1| = 1,
$$
having quasiconformal extensions across the unit circle $\mathbb S^1 =
\partial \D$ to the whole sphere $\hC = \C \cup \{\iy\}$, which satisfy $w(1) = 1$.

Equivalently, this collection is a disjunct union
$$
\wh S(1) = \bigcup_{- \pi \le \theta < \pi} S_\theta(1),
$$
where $S_\theta(1)$ consists of univalent functions $w(z) = e^{i
\theta} z + a_2 z^2 + \dots$ with quasiconformal extensions to $\hC$
satisfying $w(1) = 1$ (also completed in the indicated weak
topology). In the general case (for the limit functions of sequences of functions
with qusiconformal extension) the equality $w(1) = 1$ must be understand in terms
of the Carath\'{e}odory prime ends.

This family is closely related to the canonical class $S$ of
univalent functions $w(z$ on $\D$ normalized by $w(0) = 0, \
w^\prime(0) = 1$. Every $w \in S$ has its representative $\wh w$ in
$\wh S(1)$ (not necessarily unique) obtained by pre and post
compositions of $w$ with rotations $z \mapsto e^{i \alpha} z$ about
the origin, related by
 \be\label{1}
w_{\tau, \theta}(z) =  e^{- i \theta} w(e^{i \tau} z) \quad
\text{with} \ \ \tau = \arg z_0,
\end{equation}
where $z_0$ is a point for which $w(z_0) = e^{i \theta}$ is a common
point of the unit circle and the  boundary of domain $w(\D)$.

This implies, in particular, that the functions conformal in the closed disk
$\ov \D$ are dense in
each class $S_\theta(1)$. Such a dense subset is formed, for example, by the
images of the homotopy functions $[f]_r(z) = \fc{1}{r} f(r z)$ with real $r \in (0, 1)$ combined with rotations (1).

As is shown by Lemma 1 below, any function $w(z)$ from $\wh S(1)$ is holomorphic in the disk $\D$ (has there no pole). Their {\bf Schwarzian derivatives}
$$
S_w(z) = \left(\frac{w^{\prime\prime}(z)}{w^\prime(z)}\right)^\prime
- \frac{1}{2} \left(\frac{w^{\prime\prime}(z)}{w^\prime(z)}\right)^2,
\quad z \in \D,
$$
belong to the complex Banach space $\B$ of hyperbolically bounded
holomorphic functions  $\vp$ (more precisely, of holomorphic quadratic differentials
$\vp(z) dz^2$ on $\D$ with the norm
$$
\|\vp\|_\B = \sup_D \ld_\D^{-2}(z) |\vp(z)|,
$$
where
$\ld_\D(z) = 1/(1 - |z|^2)$.
Accordingly, $\ld_\D(z) |dz|$ is the hyperbolic metric on $\D$ of Gaussian
curvature $- 4$. This space $\B$ is dual to the space $A(\D)$ of
integrable holomorphic functions on $\D$ with $L_1$ norm.

Every $\vp \in \B$ is the Schwarzian derivative $S_w$ of a locally
univalent function $f(z)$ in the disk $\D$ determined (up to a M\"{o}bius map of the sphere $\hC$) from the nonlinear differential equation
$$
(w^{\prime\prime}/w^\prime)^\prime - (w^{\prime\prime}/w^\prime)^2/2 = \vp,
$$
or equivalently, as the ratio $w = \eta_2/\eta_1$
of two linearly independent solutions of the linear equation $2 \eta^{\prime\prime} + \vp \eta = 0$ in $\D$.

The derivatives $S_w$ of quasiconformally extendable functions $w$ from any class $S_\theta(1)$ fill in the space $\B$ a path-wise bounded domain modelling the
{\bf universal Teichm\"{u}ller space} $\T$. We identify the space $\T$ with this domain.
This allows one to consider the Taylor coefficients of functions $w$ and of their
Schwarzians as holomorphic functions on the space $\T$.

\bigskip\noindent
{\bf 1.2}.  Let $\mathcal X$ be a rotationally invariant subclass of $\wh S(1)$ obtained by completion in the topology of locally uniform convergence
in $\D$ of its dense subset $\mathcal X^0$ of quasiconformally extendable functions.
The rotational invariance means that $\mathcal X$ contains any its function
$w$  with all its pre and post rotations about the origin.

Assume also that this class satisfies the following two conditions:

$(a)${\bf openness}, which means that the corresponding collection of the Schwarzians $S_w$  of $w \in \mathcal X^0$ determines in the space $\T$  a complex Banach
submanifold (of finite or infinite dimension), which we denote by $\mathcal X_\T$;

$(b)$ {\bf variational stability}, which means that for any quasiconformal deformation $h$  of a functions $w \in \mathcal X^0$ with Beltrami coefficients $\mu_h$ supported in the complementary domain  of $w(\D)$ the composition
$h \circ w|\D$ also belongs to $\mathcal X^0$.

In fact, we shall use only a special type of such quasiconformal deformations.

We associate with such a class $\mathcal X$ the quantity
  \be\label{2}
a_2(\mathcal X) = \max \{|a_2(w)|: \ w \in \mathcal X\}
\end{equation}
and the set of rotations
 \be\label{3}
\mathcal R_\mathcal X  = \{w_{0,\tau,\theta}(z) = e^{- i \theta} w_0(e^{i \tau} z)\},
\end{equation}
where $w_0$ is one of the maximizing functions for $a_2$
on $\mathcal X$, i.e., with $|a_2(w_0)| = a_2(\mathcal X)$.

Consider on $\mathcal X$ the rotationally invariant functionals
 \begin{equation}\label{4}
J(w) = J(a_{m_1}, \dots, a_{m_s}): \ \mathcal X \to \C,
\end{equation}
which are polynomials of a distinguished set of coefficients of functions
on coefficients of  $w \in \mathcal X$.

For such classes of univalent functions and functionals, we have the following general theorem, which completely describes the extremal functions of these functionals.

\bigskip\noindent
{\bf Theorem 1}. {\it Any rotationally invariant polynomial functional (4),
whose zero set
$\mathcal Z_J = \{w \in \wh S: \ J(w) = 0\}$
is separated from the set (3), is maximized on the class $\mathcal X$ only by  functions $w_{0,\tau,\theta} \in \mathcal R_\mathcal X$.

In other words, any extremal function $w_0$ of any homogeneous (rotationally invariant) coefficient functional $J$ on a rotationally invariant and variationally
stable
class $\mathcal X$ must be simultaneously maximal for the second coefficient $a_2$ on this class, unless $J(w_0) = 0$. }

All assumptions of this theorem on the class $\mathcal X$ and the functional $J$
are essential and cannot be omitted.
This will be illustrated on examples in the last section.

\bigskip\noindent
{\bf 1.3}. The following two classes of univalent functions are of special interest.

First, let $\mathcal X$ be the canonical class $S$ of univalent functions
$w(z) = z + \sum\limits_2^\iy a_n z^n$ on $\D$ with $w(0) = 0, \ w^\prime(0) = 1$.
The classical result for this class states $|a_2| \le 2$, with equality only for
the Koebe function
 \be\label{5}
\kp_0(z) = \fc{z}{(1 - z)^2} = z + \sum\limits_2^\iy n z^n
\end{equation}
mapping the unit disk onto the complement of the ray $\{w = - t : \
1/4 \le t \le \iy\}$, and for the rotations $\kp_\theta(z) = e^{- i
\theta} \kp_0(e^{i \theta} z)$ of this function.

In this case, Theorem 1 implies an alternate and direct proof of de Branges theorem solving the Bieberbach conjecture that $|a_n|  \le n$ for all $f \in S$ (see \cite{DB}, \cite{Ha}, \cite{Kr4}), and moreover, Theorem 1 yields that the Koebe function (5)
is extremal for all coefficient functionals of type (4).

\bigskip
The second case concerns
the collections $\mathcal X(\Gamma)$ of univalent functions $w(z)$ on $\D$ compatible with the  hyperbolic Fuchsian groups $\Gamma$ acting on $\D$. This means that the maps
$$
w_\gamma = w \circ \gamma \circ w^{-1}, \quad \gamma \in \Gamma,
$$
are the Moebius transformations of the sphere $\hC$, and the group
$\Gamma^\prime = w \Gamma w^{-1}$ is a discrete subgroup the Moebius group
$\Mob (\hC)$.

As a consequence of Theorem 1, one obtains for such collections as the following result.

\bigskip\noindent
{\bf Theorem 2}. {\it (i) For every Fuchsian group $\Gamma$ and any homogeneous polynomial functional (1) on $\mathcal X(\Gamma) \subset \T$ satisfying the assumptions of Theorem 1, any its maximizing function $w_0(z)$ must simultaneously
maximize the second coefficient $a_2(w)$ on this class.

The rotations connecting the extremal functions of the functional (4) correspond
to conjugation of both groups $\Gamma$ and $w_0 \Gamma w_0^{-1}$ by elements of $\Mob(\hC)$. }

\bigskip
This theorem opens ways to investigations of extremal problems for quasiconformal deformations of Fuchsian groups.

One of the interesting open questions here is, for which groups $\Gamma$  and extremals $w_0$, the group $w_0 \Gamma w_0^{-1}$
is a totally degenerated functional group, i.e., such that its set of discontinuity
$\Omega(w_0 \Gamma w_0^{-1})$ is a simply connected domain (dense in $\hC$).

\bigskip\noindent
{\bf 1.4}. As was mentioned above, the proof of Theorem 1 involves the deep
results of Teichm\"{u}ller space theory.
The functional  $J$ is lifted from $\mathcal X$ to the
the corresponding submanifold $\mathcal X_\T$ in the Teichm\"{u}ller space
$\T_1$ of the punctured disk $\D_{*} = \{0 < |z| < 1\}$. This space is
biholomorphically equivalent to the Bers
fiber space $\mathcal F(\T)$ over the universal Teichm\"{u}ller
space $\T = \Teich (\D)$ (see \cite{Be}). This generates a holomorphic functional
$\mathcal J(\vp, t)$ on the image of $\mathcal X_\T$ in $\mathcal F(\T)$
covering $J$, with the same range domain as the initial finctional $J$. Here $\vp = S_w \in \mathcal X$ are
the Schwarzian derivatives of the initial univalent functions, while the second
variable $t$ runs over the fiber domain $w_\vp(\D)$ defined by $\vp$ .

A crucial step in the proof is to maximize $|\mathcal J(\vp, t)|$ over
$\vp$ by a fixed $t$.

\bigskip\bigskip
\centerline{\bf 2. A GLIMPSE TO TEICHM\"{U}LLER SPACES}

\bigskip
We briefly recall some needed results from Teichm\"{u}ller space
theory in order to prove our theorems; the details
can be found, for example, in \cite{Be}, \cite{GL}, \cite{Le}.

\bigskip\noindent
{\bf 2.1}.  The {\bf universal Teichm\"{u}ller space} $\T = \Teich (\D)$
is the space of quasisymmetric homeomorphisms of the unit circle
$\mathbb S^1$ factorized by M\"{o}bius maps;  all Teichm\"{u}ller
spaces have their isometric copies in $\T$.

The canonical complex Banach structure on $\T$ is defined by
factorization of the ball of the Beltrami coefficients (or complex
dilatations)
$$
\Belt(\D)_1 = \{\mu \in L_\iy(\C): \ \mu|\D^* = 0, \ \|\mu\| < 1\}
$$
vanishing on the complementary disk $\D^* = \{z \in \hC: \ |z| > 1\}$.

The coefficients $\mu_1, \mu_2 \in \Belt(\D)_1$ are called {\bf equivalent}
if the
corresponding \qc \ maps $w^{\mu_1}, w^{\mu_2}$ (solutions to the
Beltrami equation $\partial_{\ov{z}} w = \mu \partial_z w$ with $\mu
= \mu_1, \mu_2$) coincide on the unit circle $\mathbb S^1 = \partial
\D^*$ (hence, on $\ov{\D^*}$). Such $\mu$ and the corresponding maps
$w^\mu$ are called $\T$-{\it equivalent}. The equivalence classes
$[w^\mu]_\T$ are in one-to-one correspondence with the Schwarzian
derivatives $S_{w^\mu}(z), \ z \in \D^*$, which belong to the space
$\B = \B(\D^*)$ of hyperbolically bounded holomorphic functions on
the disk $\D^*$  with norm
$$
\|\vp\|_\B = \sup_{D^*} (|z|^2 - 1)^2 |\vp(z)|.
$$
Note that $\vp(z) = O(|z| ^{-4})$ as $z \to \infty$.

This space is dual to the Bergman space $A_1(\D^*)$, a subspace
of $L_1(\D^*)$ formed by integrable holomorphic functions (quadratic
differentials $\vp(z) dz^2$ on $D$), since every linear functional
$l(\vp)$ on $A_1(D)$ is represented in the form
$$
l(\vp) = \langle \psi, \vp \rangle_{\D^*} = \iint\limits_{\D^*} \
(|z|^2 - 1)^2 \ov{\psi(z)} \vp(z) dx dy
$$
with a uniquely determined $\psi \in \B(\D^*)$.

The Schwarzians $S_{w^\mu}(z)$ with $\mu \in \Belt(\D)_1$ range over
a bounded domain in the space $\B = \B(\D^*)$. This domain models
the space $\T$. It lies in the ball $\{\|\vp\|_\B < 6\}$ and
contains the ball $\{\|\vp\|_\B < 2\}$. In this model, the
Teichm\"{u}ller spaces of all hyperbolic Riemann surfaces are
contained in $\T$ as its complex submanifolds.

The factorizing projection
$$
\phi_\T(\mu) = S_{w^\mu}: \ \Belt(\D)_1 \to \T
$$
is a holomorphic map from $L_\iy(\D)$ to $\B$. This map is a split
submersion, which means that $\phi_\T$ has local holomorphic
sections (see, e.g., \cite{EKK}, \cite{GL}).

Both equations $S_w = \vp$ and $\partial_{\ov z} w = \mu
\partial_z w$ (on $\D^*$ and $\D$, respectively) determine their
solutions up to a M\"{o}bius transformation of $\hC$.
So appropriate normalization of solution $w^\mu(z)$ (for example,
fixing the points $1, i, - 1$ or other three points on the unit circle),
provides uniqueness of solution
of either equation, and moreover, then the values $w^\mu(z_0)$ at any
point $z_0 \in \C \setminus \{1, i, -1\}$ and the Taylor coefficients
$b_1, b_2, \dots$ of $w^\mu \in \Sigma_\theta$ depend holomorphically
on $\mu \in \Belt(\D)_1$ and on $S_{w^\mu} \in \T$.
Later we shall use another normalization which also insures the needed
uniqueness and holomorphy.

\bigskip\noindent
{\bf 2.2}. The points of Teichm\"{u}ller space $\T_1  =  \Teich(\D_{*})$
of the punctured disk $\D_{*} =  \{0 < |z| < 1\}$
are the classes $[\mu]_{\T_1}$ of $\T_1$-{\bf equivalent} \Be \
coefficients $\mu \in \Belt(\D)_1$ so that the corresponding \qc \
automorphisms $w^\mu$ of the unit disk coincide on both boundary
components (unit circle $\mathbb S^1 = \{|z| =1\}$ and the puncture
$z = 0$) and are homotopic on $\D \setminus \{0\}$. This space can
be endowed with a canonical complex structure of a complex Banach
manifold and embedded into $\T$ using uniformization.

Namely, the disk $\D_{*}$ is conformally equivalent to the factor
$\D/\G$, where $\G$ is a cyclic parabolic Fuchsian group acting
discontinuously on $\D$ and $\D^*$. The functions $\mu \in
L_\iy(\D)$ are lifted to $\D$ as the \Be \ $(-1, 1)$-measurable
forms  $\wt \mu d\ov{z}/dz$ in $\D$ with respect to $\G$, i.e., via
$(\wt \mu \circ \g) \ov{\g^\prime}/\g^\prime = \wt \mu, \ \g \in
\G$, forming the Banach space $L_\iy(\D, \G)$.

We extend these $\wt \mu$ by zero to $\D^*$ and consider the unit
ball $\Belt(\D, \G)_1$ of $L_\iy(\D, \G)$. Then the corresponding
Schwarzians $S_{w^{\wt \mu}|\D^*}$ belong to $\T$. Moreover, $\T_1$
is canonically isomorphic to the subspace $\T(\G) = \T \cap \B(\G)$,
where $\B(\G)$ consists of elements $\vp \in \B$ satisfying $(\vp
\circ \g) (\g^\prime)^2 = \vp$ in $\D^*$ for all $\g \in \G$.

Due to the Bers isomorphism theorem, the space $\T_1$ is
biholomorphically isomorphic to the Bers fiber space
$$
\mathcal F(\T) = \{(\phi_\T(\mu), z) \in \T \times \C: \ \mu \in
\Belt(\D)_1, \ z \in w^\mu(\D)\}
$$
over the universal space $\T$ with holomorphic projection $\pi(\psi,
z) = \psi$ (see \cite{Be}).

This fiber space is a bounded hyperbolic domain in $\B \times \C$
and represents the collection of domains $D_\mu = w^\mu(\D)$ as a
holomorphic family over the space $\T$. For every $z \in \D$,  its
orbit $w^\mu(z)$ in $\T_1$ is a holomorphic curve over $\T$.

The indicated isomorphism between $\T_1$ and $\mathcal F(\T)$ is
induced by the inclusion map \linebreak $j: \ \D_{*} \hookrightarrow
\D$ forgetting the puncture at the origin via
 \be\label{6}
\mu \mapsto (S_{w^{\mu_1}}, w^{\mu_1}(0)) \quad \text{with} \ \
\mu_1 = j_{*} \mu := (\mu \circ j_0) \ov{j_0^\prime}/j_0^\prime,
\end{equation}
where $j_0$ is the lift of $j$ to $\D$.

By Koebe's one-quarter theorem, for any univalent function $W(z) =
z + b_0 + b_1 z^{-1} + \dots$ in $\D^*$, the boundary of domain $W(D*)$
is located in the disk $\{|w - b_0| \le 2\}$. If $W(z) \ne 0$ in $\D^*$,
its inversion $w(z) = z + a_2 z^2 + \dots$ is univalent in $\D$, and
$b_0 = - a_2$ satisfies $|b_| \le 2$.
Using the maps $W$ with
quasiconformal extensions, one gets by the Bers theorem that the
indicated domains $D_\mu$ are filled by the admissible values of $W^\mu(0)$;
all these domains are located in the disk $\{|W| \le 4\}$.

In the line with our goals, we slightly modified the Bers
construction, applying quasiconformal maps $F^\mu$ of $\D_{*}$
admitting conformal extension to $\D^*$ (and accordingly using  the
Beltrami coefficients $\mu$ supported in the disk) (cf. \cite{Kr2}).
These changes are not essential and do not affect the underlying
features of the Bers isomorphism (giving the same space up to a
biholomorphic isomorphism).

The Bers theorem is valid for Teichm\"{u}ller spaces $\T(X_0
\setminus \{x_0\})$ of all punctured hyperbolic Riemann surfaces
$X_0 \setminus \{x_0\}$ and implies that $\T(X_0 \setminus \{x_0\})$
is biholomorphically isomorphic to the Bers fiber space $ \mathcal
F(\T(X_0))$ over $\T(X_0)$.

Note also that every Teichm\"{u}ller space $\T(X)$ is a complete metric space with intrinsic Teichm\"{u}ller metric $\tau_\T(\cdot, \cdot)$ defined by quasiconformal maps. By the
Royden-Gardiner theorem, this metric is equal the hyperbolic Kobayashi metric
$d_\T(\cdot, \cdot)$
determined by the complex structure on this space (see, e.g., \cite{EKK}, \cite{GL}, \cite{Ro}). In other words, the Kobayashi-Teichm\"{u}ller metric is the maximal invariant metric on $\T(X)$.

\bigskip\bigskip
\centerline{\bf 3. UNDERLYING LEMMA}

\bigskip
We start with the following lemma, which ensures the existence of univalent functions in the disk with quasiconformal extension satisfying the prescribed normalization of classes $S_\theta(1)$
and some other conditions. It concerns the solutions of the Beltrami equation
$\partial_{\ov z} w = \mu(z) \partial_z w$ on $\C$ with
coefficients $\mu$ supported in the disk $\D^*$, i.e., from the ball
$$
\Belt(\D^*)_1 = \{\mu \in L_\iy(\C): \ \mu|\D = 0, \ \|\mu\| < 1\}
$$
(and hence the solutions of the corresponding Schwarzian equation $S_w(z) = \vp$ in $\D$ with given $\vp \in \B$).

\bigskip\noindent
{\bf Lemma 1}. {\it For any $\mu \in \Belt(\D^*)_1$ and any $\theta \in [0, 2 \pi)$, there exists a unique homeomorphic solution $w = w^\mu(z)$ of the equation $\partial_{\ov z} w = \mu(z) \partial_z w$ on $\hC$ such that
 \be\label{7}
w(0) = 0, \quad w^\prime(0) = e^{i \theta}, \quad w(1) = 1.
\end{equation}
This solution is holomorphic on the unit disk $\D$, and hence, $w(z_0) = \infty$ at some point $z_0$ with $|z_0| \ge 1$ (so $w(z)$ does not have a pole in $\D$).  }

\bigskip\noindent
{\bf Proof}. Let us first consider the coefficients $\mu$ vanishing in a broader
disk $\D_r = \{|z| < r\}, \ r > 1$, so that $w^\mu$ is conformal on $\D_r \Supset \ov{\D}$, and  assume that $\mu \neq \mathbf 0$ (the origin of $\Belt(\D^*)_1$).

Fix $a \in (1, r) $ close to $1$ and $\theta \in [0, 2 \pi]$; then $1/a \in \D$.

The generalized Riemann mapping theorem for the Beltrami equation
$\partial_{\ov z} w = \mu(z) \partial_z w$  on $\hC$
implies a homeomorphic solution $\wh w$ to this equation satisfying
 \be\label{8}
\wh w(- 1/a) = - 1/a, \quad \wh w^\prime(- 1/a) = e^{i \theta} \ , \quad \wh w(\infty) = \infty.
\end{equation}
Its composition with the M\"{o}bius map
$$
\gamma_a(z) = (1 - a z)/(z - a)
$$
preserving either from disks $\D$ and $\D^*$ has the Beltrami coefficient
$$
\gamma_{a,*} \mu :=  \mu_{\wh w \circ \gamma_a}(z) = \mu
\circ \gamma_a(z) \gamma_a^\prime(z)/\ov{\gamma_a^\prime(z)}
$$
and also is conformal in the disk $\D_r$ and holomorphic in $\D$.

Since, by the classical Schwarz lemma, for any holomorphic map
$g: \D \to \D$ and any point $z_0 \in \D$,
$$
|g^\prime(z_0)| \le (1- |g(z_0)|^2)/(1 - |z_0|^2)
$$
with equality only for appropriate M\"{o}bius automorphism of $\D$, the above  normalization (8) and the assumption on $\mu$ yield that for any map $\wh w(z)$ the image $\wh w(\D)$ does not cover $\D$, and thus either $\wh w(\D)$ is a proper subdomain of $\D$ or it also contains  the points $z$ with $|z| > 1$ outer for $\D$.

Passing if needed to suitable rotated map $e^{-i \a} \wh w(e^{i \a} z)$, one obtains
that this domain $\wh w(\D)$ does not contain simultaneously both distinguished points
$a$ and $1/a$ (at least sufficiently close to $1$).

Now consider the map
$$
w_{a,a}(z) = \gamma_a^{-1} \circ \wh w \circ \gamma_a(z),
$$
having the same Beltrami coefficient $\gamma_{a,*} \mu$. Since
$$
\gamma_a(\infty) = - a, \quad \gamma_a(a) = \infty, \quad \gamma_a(0) = - 1/a,
\quad \gamma_a(1/a) = 0
$$
(and accordingly, $\gamma_a^{-1}(\infty) = a, \ \gamma_a^{-1}(- a) = \infty, \ \gamma_a^{-1}(0) = 1/a, \ \gamma_a^{-1}(- 1/a) = 0$), the map $w_{a,a}$  satisfies
 \be\label{9}
w_{a,a}(0) = 0, \quad w_{a,a}^\prime(0) = \wh w^\prime(- 1/a) = e^{i \theta}, \quad
w_{a,a} (a) = \gamma_a^{-1} \circ w \circ \gamma_a(a) = a.
\end{equation}

If one starts with Beltrami coefficient $\gamma_{a,*}^{-1} \mu$, taking its map
$w$ normalized by (7), then the final map $w_{a,a}$ has the initial Beltrami coefficient
$\mu$ and satisfies (8) for all $a \in (1, r)$ sufficiently close to $1$.

In view of our assumptions on $\wh w$, the point $\wh w^{-1}(a)$ cannot lie in
the unit disk $\D$;
therefore the function $w_{a,a}$ is holomorphic in this disk.

\bigskip
Now we investigate the limit process as $a \to 1$.
Any from the constructed maps $w_{a,a}$ is represented as a composition of a fixed solution $\wh w$ to the equation
$\partial_{\ov z} w = \mu(z) \partial_z w$ subject to (10) and some  M\"{o}bius maps $\wh \gamma_a$. The first two conditions in (9)
imply that the restrictions of these $\wh \gamma_a$ to $\wh w(\D_r)$ form a
(sequentially) compact set of $\wh \gamma_a$ in the topology of convergence in the spherical metric on $\hC$. Letting $a \to 1$, one
obtains in the limit the map $\wh \gamma_1(z) = \lim\limits_{a\to 1} \wh \gamma_a(z)$,
which also is a non-degenerate (nonconstant) M\"{o}bius map. Accordingly,
$$
\lim\limits_{a\to 1} w_{a,a}(z) = \wh \gamma_1 \circ \wh w(z) =: \wh w_1(z),
$$
and this map satisfies (9) with $a = 1$, which is equivalent to (7).

Note that the relations (9) do not depend on $r$ and that the normalization (8) also is valid for the inverse rotation $e^{i \a} \wh w_1(e^{-i \a} z)$ of the limit function.

This implies the assertion of Lemma 2 for all Beltrami coefficients $\mu \ne \mathbf 0$ supported in the disk $\D_r^* = \{|z| > r\}$ with $r > 1$.

\bigskip
To extend the obtained result to arbitrary $\mu \in \Belt(\D^*)_1 \ (\mu \ne \mathbf 0$), we pass to coefficients $\mu_r(z) = \mu(z)$ for $|z| > r$ and extended by zero to $\D_r$. The compactness properties of the $k$-quasiconformal families (i.e., with
$\|\mu\|_\iy \le k < 1$)
imply the convergence of maps $w^{\mu_r}(z)$ normalized by (7) to $w^\mu$(z) as $r \to 1$
in the spherical metric on $\hC$ (and hence everywhere on $\hC$).

The remained case $\mu(z) \equiv 0$ omitted above follows in the limit as $\mu \to \mathbf 0$ (or even $\mu(z) \to 0$ almost everywhere in $\D^*$). Then the map $w(z)$ satisfying (7) must be an
elliptic fractional linear transformation with fixed points $0$ and $1$; hence,
$$
\frac{w - 1}{w} = e^{i\theta} \ \frac{z - 1}{z},
$$
which implies
  \be\label{10}
w = \frac{e^{-i \theta} z}{(e^{-i \theta} - 1) z + 1}.
\end{equation}
The simple direct calculations yield that $w(z_0) = \infty$ only at a point $z_0$
with $|z_0| \ge 1$ (which also follows from the above).
The proof of Lemma 2 is completed.

This lemma plays an important role in our further considerations. We shall consider
the univalent function $w(z)$ in the disk $\D$ normalized by (9) and their rotations
(1) with $\tau, \theta \in [0, 2 \pi]$; all these rotations also are holomorphic univalent in this disk, so $w_{\tau,\theta}(z_0) = \infty$ only at some point $z_0 \in \ov{\D^*}$.

\bigskip\bigskip
\centerline{\bf 4. PROOF OF THEOREM 1}

\bigskip
We accomplish the proof of this theorem in four stages.

\bigskip\noindent
{\bf Step 1: Lifting the functional $J(f)$ onto the universal Teichm\"{u}ller space}.
The prescribed normalizing conditions
$w(0) = 0, \ w^\prime(0) = e^{i \theta}, \ w(1) = 1$
are compatible with existence and uniqueness of the corresponding conformal and  quasiconformal maps and the Teichm\"{u}ller space theory, ensure holomorphy of their Taylor coefficients, etc.
Actually we deal with the classical model of Teichm\"{u}ller spaces via
domains in the Banach spaces of Schwarzian dervatives $S_w$ in $\mathbb D$
(or in the disk $\mathbb D^*$) of univalent holomorphic functions normalized
either by fixing three boundary points on the unit circle $S^1$ or via $w(0) = 0, \ w^\prime(0) = 1, \ w(z_0) = z_0$, where $z_0$ on $\mathbf S^1$.
(Often the disk is replaced by the half-plane.)

Modelling the universal Teichm\"{u}ller space $\mathbf T$ by the Schwarzians
$S_w = \vp$ of functions $w(z)$ from $S_\theta(1)$, we have that its base point
$\varphi = \mathbf 0$  corresponds to the function (10) (which equals the identity map for $\theta = 0$).

It is more convenient technically to deal with univalent functions in
the complementary disk $\D^*$. In view of Lemma 1, we can model the space $\T$
using the inverted functions $W(z) = 1/w(1/z)$ for $w \in \wh S(1)$.

These functions form the corresponding classes $\Sigma_\theta(1)$ of nonvanishing
univalent functions on the disk $\D^*$ with expansions
$$
W(z) =  e^{- i \theta} z + b_0 + b_1 z^{-1} + b_2 z^{-2} + \dots,
\quad  W(1) = 1,
$$
and $\wh \Sigma(1) = \bigcup_\theta \Sigma_\theta(1)$.
It is more convenient technically to deal with these functions.

Simple computations yield that the coefficients $a_n$ of $f \in
S_\theta(1)$ and the corresponding coefficients $b_j$ of $W(z) =
1/f(1/z) \in \Sigma_\theta(1)$ are related by
$$
b_0 + e^{2i \theta} a_2 = 0, \quad b_n + \sum \limits_{j=1}^{n}
\epsilon_{n,j}  b_{n-j} a_{j+1} + \epsilon_{n+2,0} a_{n+2} = 0,
\quad n = 1, 2, ... \ ,
$$
where $\epsilon_{n,j}$ are the entire powers of $e^{i \theta}$. This
successively implies the representations of $a_n$ by $b_j$ via
 \be\label{11}
a_n = (- 1)^{n-1} \epsilon_{n-1,0}  b_0^{n-1} - (- 1)^{n-1} (n - 2)
\epsilon_{1,n-3} b_1 b_0^{n-3} + \text{lower terms with respect to}
\ b_0.
\end{equation}

This transforms the initial functional (4) into a coefficient functional $\wt J(W)$
on $\Sigma_\theta(1)$ depending on the corresponding coefficients $b_j$.

Note that the coefficients $\a_n$ of Schwarzians
$S_w(z) = \sum_0^\infty \a_n z^n$
are represented as polynomials of $n + 2$ initial coefficients of $w
\in S_\theta(1)$ and, in view of (11), as polynomials of $n + 1$
initial coefficients of the corresponding $W \in \Sigma_\theta(1)$,
provided that $\theta$ is given and fixed and the number $e^{i \theta}$
is considered to be a constant (vice versa, the coefficients $a_n$ and $b_j$ are uniquely determined by $\a_n$ by solving the Schwarzian differential equation
$S_w = \vp$ or from the equation $S_W = \vp(1/z)$ and (11).

We will deal with our polynomial functionals $J(w)$ and $\wt J(W)$
only on a fixed class $S_\theta(1)$ or $\Sigma_\theta(1)$.

\bigskip
Holomorphic dependence of normalized quasiconformal maps on complex
parameters (first established by Ahlfors and Bers in \cite{AB} for maps with
three fixed points on $\hC$)
is an underlying fact for the Teichm\"{u}ller space theory and for many
other applications.

Another somewhat equivalent proof of holomorphy involves the variational
technique for quasiconformal maps. For the maps $w$ from $S_\theta(1)$,
this holomorphy is a consequence of the following lemma  from \cite{Kr1}, Ch. 5,
combined with appropriate M\"{o}bius maps.

\bigskip\noindent
{\bf Lemma 2}. {\it Let $w(z)$ be a  quasiconformal map of the plane
$\hC$ with Beltrami coefficient $\mu(z)$ which satisfies
$\|\mu\|_\iy < \ve_0 < 1$ and vanishes in the disk $\{|z| < r\}$.
Suppose that $w(0) = 0, \ w^\prime(0) = 1$, and $w(1) = 1$. Then,
for sufficiently small $\ve_0$ and for $|z| \le R < r_0(\ve_0, r)$
we have the variational formula
$$
w(z) = z - \frac{z^2 (z - 1)}{\pi} \iint\limits_{|\zeta|>r} \
\frac{\mu(\zeta) d\xi d \eta}{\zeta^2(\zeta -1)(\zeta- z)} +\Omega_\mu(z),
$$
where $\zeta = \xi + i \eta; \
\max_{|z|\le R} |\Omega_\mu(z) \le C(\ve_0, r, R) \|\mu\|_\iy^2; \
r_0(\ve_0, r)$ is a well defined function of $\ve_0$ and $r$ such that
$\lim_{\ve_0\to 0} \ r_0(\ve_0, r) = \infty$, and the constant $C(\ve_0, r, R)$
depends only on $\ve_0, \ r$ and $R$.   }

\bigskip\noindent
{\bf Step 2: Lifting to covering space $\T_1$ and estimating the
restricted plurisubharmonic functional}.
Our next step is to lift both polynomial functionals $J(w)$ and
$\wh J(W)$ onto the Teichm\"{u}ller space $\T_1$ which covers $\T$.

Letting
 \be\label{12}
\wh J(\mu) = \wt J(W^\mu),
\end{equation}
we lift these functionals  from the sets $S_\theta(1)$ and
$\Sigma_\theta(1)$ onto the ball $\Belt(\D)_1$. Then, under the
indicated $\T_1$-equivalence, i.e., by the quotient map
$$
\phi_{\T_1}: \ \Belt(\D)_1 \to \T_1, \quad \mu \to [\mu]_{\T_1},
$$
the functional $\wt J(W^\mu)$ is pushed down to a bounded holomorphic
functional $\mathcal J$ on the space $\T_1$ with the same range domain.

Equivalently, one can apply the quotient map $\Belt(\D)_1 \to \T$
(i.e., $\T$-equivalence) and compose  the descended functional on
$\T$ with the natural holomorphic map $\iota_1: \ \T_1 \to \T$
generated by the inclusion $\D_{*} \hookrightarrow \D$ forgetting
the puncture. Note that since the coefficients $b_0, \ b_1, \dots$
of $W^\mu \in \Sigma_\theta$   are uniquely determined by its
Schwarzian $S_{W^\mu}$, the values of $\mathcal J$ in the points
$X_1, \ X_2 \in \T_1$ with $\iota_1(X_1) = \iota_1(X_2)$ are equal.

Now, using the Bers isomorphism theorem, we regard the points of the
space $\T_1$ as the pairs $X_{W^\mu} = (S_{W^\mu}, W^\mu(0))$, where
$\mu \in \Belt(\D)_1$ obey $\T_1$-equivalence (hence, also
$\T$-equivalence). Denote (for simplicity of notations) the
composition of $\mathcal J$ with biholomorphism $\T_1 \cong \mathcal
F(\T)$ again by $\mathcal J$. In view of (5) and (13), it is
presented on the fiber space $\mathcal F(\T)$ by
 \be\label{13}
\mathcal J(X_{W^\mu}) = \mathcal J(S_{W^\mu}, \ t), \quad t = W^\mu(0).
\end{equation}
This yields a logarithmically plurisubharmonic functional
$|\mathcal J(S_{W^\mu}, t)|$ on $\mathcal F(\T)$.

Note that since the coefficients $b_0, \ b_1, \dots$ of $W^\mu \in
\Sigma_\theta$   are uniquely determined by its Schwarzian
$S_{W^\mu}$, the values of $\mathcal J$ in the points $X_1, \ X_2
\in \T_1$ with $\iota_1(X_1) = \iota_1(X_2)$ are equal.

\bigskip
We have to estimate a smaller plurisubharmonic functional arising
after restriction of $\mathcal J(S_{W^\mu}, \ t)$ to $S_W \in \mathcal X_\T$,
i.e., the restriction of functional (13) onto the corresponding set of
pairs $(S_{W^\mu}, W^\mu(0))$
consisting of $S_{W^\mu} \in \mathcal X_\T$ and of the values $W^\mu(0)$ filling some subdomain $D_{\mathcal X,\theta}$.

Since our functionals are polynomials, they are defined for all $S_W \in \T$ and $t$ from some domain $D_\theta$ containing $D_{\mathcal X,\theta}$.
We define on $D_\theta$ the function
 \be\label{14}
u_\theta(t) = \sup_{W^\mu} |\mathcal J(S_{W^\mu}, t)|,
\end{equation}
where the supremum is taken over all $S_{W^\mu} \in \T$ admissible for a given $t =
W^\mu(0) \in D_\theta$.

The following basic lemma from \cite{Kr5} provides that this function inherits subharmonicity of $\mathcal J$.

\bigskip\noindent
{\bf Lemma 3}. {\it The function $u_\theta(t)$ is subharmonic in the domain $D_\theta$. }

\bigskip
The proof of this lemma in \cite{Kr5} is complicated. It involves a weak approximation of the underlying space $\T$ (and simultaneously of the space $\T_1$) by finite dimensional Teichm\"{u}ller spaces of the punctured spheres in the topology of locally uniform convergence on $\C$ and using the increasing unions of the quotient spaces
 \be\label{15}
\mathcal T_s = \bigcup_{j=1}^s \ \wh \Sigma_{\theta_j}^0/\thicksim \
= \bigcup_{j=1}^s \{(S_{W_{\theta_j}}, W_\theta^\mu(0)) \} \ \simeq
\T_1 \cup \dots \cup \T_1,
\end{equation}
where $\theta_j$ run over a dense subset $\Theta \subset [-\pi, \pi]$,
the equivalence relation $\thicksim$ means $\T_1$-equivalence
on a dense subset $\wh \Sigma^0(1)$ in the union $\wh \Sigma(1)$
formed by univalent functions $W_{\theta_j}(z) = e^{-i \theta_j} z +
b_0 + b_1 z^{-2} + \dots$ on $\D^*$ with quasiconformal extension to
$\hC$ satisfying $W_{\theta_j}(1) = 1$, and
$$
\mathbf W_\theta^\mu(0) := (W_{\theta_1}^{\mu_1}(0), \dots ,
W_{\theta_s}^{\mu_s}(0)).
$$
The Beltrami coefficients  $\mu_j \in \Belt(\D)_1$ are chosen here
independently. The corresponding collection
$$
\beta = (\beta_1, \dots, \beta_s)
$$
of the Bers isomorphisms
$$
\beta_j: \ \{(S_{W_{\theta_j}}, W_{\theta_j}^{\mu_j}(0))\} \to
\mathcal F(\T)
$$
determines a holomorphic surjection of the space $\mathcal T_s$
onto $\mathcal F(\T)$. The function (14) is determined by
$$
u(t) = \sup_{\Theta} u_{\theta_s}(t),
$$
where $u_{\theta_s}$ is obtained by maximization of type (14) over $\mathcal T_s$.

Since, by assumption, the image $\mathcal X_\T$ of the class $\mathcal X$ in $\T$ is
a complex manifold, the restriction of the function $|\mathcal J(S_{W^\mu}, t)|$
to this manifold and to corresponding values of $t = W^\mu(0)$ also is
plurisubharmonic. The arguments from \cite{Kr5} are straightforwardly extended to this restriction, giving in a similar way, that also the corresponding maximal function
$$
u_\theta(t) = \sup_{W^\mu} |\mathcal J(S_{W^\mu}, t)|
$$
is subharmonic on the domain $D_{\mathcal X, \theta}$.

\bigskip\noindent
{\bf Step 3: Range domain of $W^\mu(0)$}.
The next step in maximization of the function $u_\theta$ (and thereby of the
functional $\mathcal J$) is to establish the value domain of $W^\mu(0)$ for $W^\mu$ running over $\mathcal X_\T$. This requires the corresponding covering estimate.

Let $G$ be a domain in a complex Banach space $X = \{\mathbf x\}$
and $\chi$ be a holomorphic map from $G$ into the universal
Teichm\"{u}ller space $\T$ modeled as a bounded subdomain of $\B$.
Consider in the unit disk the corresponding Schwarz differential
equations
$S_w(z) = \chi(\x)$
and pick their holomorphic univalent solutions $w(z)$ in $\D$ satisfying $w(0) = 0, \
w^\prime(0) = 1$ (hence $w(z) = z  + \sum_2^\infty a_n z^n$).
Put
 \be\label{16}
|a_2^0| = \sup \{|a_2|: \ S_w \in \chi(G)\},
\end{equation}
and let $w_0(z) = z + a_2^0 z^2 + \dots$ be one of the maximizing
functions.

\bigskip\noindent
{\bf Lemma 4}. {\it (a) For every indicated solution
$w(z) = z + a_2 z^2 + \dots$ of the Schwarz differential equation,
the image domain $w(\D)$ covers entirely the disk
$D_{1/(2 |a_2^0|)} = \{|w| < 1/(2 |a_2^0|)\}$.

The radius value $1/(2 |a_2^0|)$ is sharp for this collection of
functions, and the circle $\{|w| = 1/(2 |a_2^0|)$ contains points
not belonging to $w(\D)$ if and only if $|a_2| = |a_2^0|$
(i.e., when $w$ is one of the maximizing functions).

(b) The inverted functions
$$
W(\zeta) = 1/w(1/\zeta) = \zeta - a_2^0
+ b_1 \zeta^{-1} + b_2 \zeta^{-2} + \dots
$$
map the disk $\D^*$ onto a domain whose boundary is entirely
contained in the disk} $\{|W + a_2^0| \le |a_2^0|\}$.

\bigskip\noindent
{\bf The proof} follows the classical lines of Koebe's $1/4$
theorem (cf. \cite{Go}).

{\it (a)} Suppose that the point $w = c$ does not belong to the
image of $\D$ under the map $w(z)$ defined above. Then $c \ne 0$,
and the function
$$
w_1(z) = c w(z)/(c - w(z)) = z + (a_2 + 1/c) z^2 + \dots
$$
also belongs to this class, and hence by (16), $|a_2 +1/c| \le
|a_2^0|$, which implies
$$
|c| \ge 1/(2 |a_2^0|).
$$
The equality holds only when
$$
|a_2 + 1/c| = |1/c| - |a_2| = |a_2^0| \quad \text{and} \ \ |a_2| =
|a_2^0|.
$$

{\it(b)} If a point $\zeta = c$ does not belong to the image
$W(\D^*)$, then the function
$$
W_1(z) = 1/[W(1/z) - c] = z + (c + a_2) z^2 + \dots
$$
is holomorphic and univalent in the disk $\D$, and therefore,  $|c +
a_2| \le |a_2^0|$. The lemma follows.

This lemma implies that the boundary of the range domain of
$W^\mu(0)$ is contained in the disk
 \be\label{17}
\D_{2|a_2^0|} = \{W: \ |W| < 2 |a_2^0|\}
\end{equation}
and the boundary of this domain touches from inside the circle
$\{|W| = 2 |a_2^0|\}$ at the points corresponding to extremal
functions maximizing $|a_2|$ on the closure of the manifold $\mathcal X_\T$.

\bigskip\noindent
{\bf Step 5. Finishing the proof}.
To complete the proof of the theorem, we have to apply some special variations
of univalent functions with quasiconformal extension given by the following
lemma, which is a special case of more general results from \cite{Kr1}.
Here we essentially use the variational stability of $\mathcal X$.

\bigskip\noindent
{\bf Lemma 5}. {\it Let $D$ be a simply connected domain on the Riemann sphere
$\hC$. Assume that there are a set $E$ of positive two-dimensional
Lebesgue measure and a finite number of points
 $z_1, z_2, ..., z_m$ distinguished in $D$. Let
$\a_1, \a_2, ..., \a_m$ be non-negative integers assigned to $z_1,
z_2, ..., z_m$, respectively, so that $\a_j = 0$ if $z_j \in E$.

Then, for a sufficiently small $\ve_0 > 0$ and $\varepsilon \in (0,
\varepsilon_0)$, and for any given collection of numbers $w_{sj}, s
= 0, 1, ..., \a_j, \ j = 1,2, ..., m$ which satisfy the conditions
$w_{0j} \in D$, \
$$
|w_{0j} - z_j| \le \ve, \ \ |w_{1j} - 1| \le \ve, \ \ |w_{sj}| \le
\ve \ (s = 0, 1, \dots   a_j, \ j = 1, ..., m),
$$
there exists a quasiconformal automorphism $h$ of $D$ which is conformal on $D
\setminus E$ and satisfies
$$
h^{(s)}(z_j) = w_{sj} \quad \text{for all} \ s =0, 1, ..., \a_j, \ j
= 1, ..., m.
$$
Moreover, the Beltrami coefficient $\mu_h(z) = \partial_{\bar z}
h/\partial_z h$ of $h$ on $E$ satisfies $\| \mu_h \|_\iy \leq M
\ve$. The constants $\ve_0$ and $M$ depend only upon the sets $D, E$
and the vectors $(z_1, ..., z_m)$ and $(\a_1, ..., \a_m)$.

If the boundary $\partial D$ is Jordan or is $C^{l + \a}$-smooth,
where $0 < \a < 1$ and $l \geq 1$, we can also take $z_j \in
\partial D$ with $\a_j = 0$ or $\a_j \leq l$, respectively.   }

\bigskip
We apply this lemma to quasiconformally extendable functions from $\mathcal X$, which are dense in this subclass (in topology of locally uniform convergence on $\D$).

Let $w_0$ be an extremal of a given functional $J$ on $\mathcal X$. Pass to
the function
$$
w_{0r}(z) = \fc{1}{r} w_0(r z) = r e^{i \theta} z + r^2 a_2^0 z^2 + \dots
$$
with $r$ close to $1$, and to its image $\wt w_{0r}$ in $\wh S(1)$, using the corresponding rotations of type (1). This image is univalent and holomorphic on the closed disk $\ov{\D}$ and satisfy also the third normalization condition $w(1) = 1$.

Varying appropriately the coefficients $a_2, a_{m_1}, \dots, a_{m_s}$ of $\wt w_{0r}$ by Lemma 5 with quasiconformal variation $h$ conformal on $\wt w_{0r}(\D)$,
one derives that the maximal function
  \be\label{18}
u(t) = \sup_\theta u_\theta(t) = \sup \big\{|\mathcal J(S_{W^\mu}, t)|: \
S(W^\mu) \in \bigcup_s \mathcal T_s \big\}
\end{equation}
(where $u_\theta$ and $ \mathcal T_s$ are determined by (14) and (15)) must be  positive on any circle $\{|t| = r < 2 |a_2^0|\}$, and hence on the whole disk (17).
This function also is subharmonic and circularly invariant on this disk.

Therefore, the maximal value of the function (18), which coincides with
$\max |J(w)|$ on $\wh S(1)$, is attained on the boundary circle $\{|t| =  2 |a_2^0|\}$ and equals $u(1)$. This is equivalent to assertion of Theorem 1.

\bigskip\bigskip
\centerline{\bf 5. ADDITIONAL REMARKS AND COUNTEREXAMPLES}

\bigskip\noindent
{\bf 1}. As was already  mentioned above, the assumptions of Theorem 1 on classes of functions and on functionals are essential and cannot be omitted.
We illustrate this on two examples.

\bigskip\noindent
{\bf Example 1}. Consider the class $S_q(\infty)$ of univalent functions
$w(z) = z + a_2 z^2 + \dots $ in $\D$ with $q$-quasiconformal extensions to $\hC \ (q < 1)$, with fixed point at infinity.

As is well known, the maximizing function for $a_2$ in $S_q(\infty)$  is
$$
w_{1,t}(z) = z/(1 - t z)^2, \quad |t| = q, 
$$ 
and $|a_2| \le 2q$. 
But for any $n > 2$ and for all $q \le 1/(n^2 + 1)$, the sharp estimate for the functional $J_n(w) = a_n$ on $f \in S_q(\iy)$ is given by 
$$
|a_n| \le 2q/(n - 1),
$$
and the equality is attained only on the functions
$$
w_{n-1,t}(z) = w_{1,t}(z^{n-1})^{1/(n-1)} = z + \fc{2 t}{n - 1} z^n
+ \dots 
$$
with $|t| = q$ (see, e.g., \cite{Kr4}). 

The reason is that this class satisfies the assumption $(a)$, but is not variationally stable, because
the generic variations of Lemma 5 provide $q^\prime$-quasiconformal maps with $q^\prime > q$.

All this is in accordance with the well known fact that {\it in any class of univalent functions with $q$-quasiconformal extension, no function can be simultaneously extremal for different holomorphic functionals  unless these functionals have equal
$1$-jets at the origin}; in particular, this holds for functionals $J_n(f) = a_n$
For details and related results we refer, for example, to \cite{Gu}, \cite{Kr3},
\cite{Kr4}, \cite{KK}, \cite{Ku}, \cite{Sc}.

The underlying feature, which causes this fact, is that for $q < 1$ these functionals $J_n$ determine the extremal disks for the Carath\'{e}odory metric on the universal Teichm\"{u}ller space $\T$.

\bigskip\noindent
{\bf Example 2.}
The class $S(M)$ formed by bounded functions $f(z) = z + a_2 z^2 + \dots \in S$
with $|f(z)| < M$ in $\D \ (M > 1)$ also is not variationally stable, because
the variations given by Lemma 5 (with the sets $E$ of quasiconformality located
outside of the disk $\D_M = \{|z| < M|\}$) generically increase the $\sup$ norm of varied functions.

Thus Theorem 1 cannot be applied also to this class. Note that this is in accordance with the fact that the known coefficient estimates for $S(M)$ (see, e.g. \cite{Pr}) are of completely different nature then Theorem 1.

\bigskip\noindent
{\bf 2}. The extension of Theorem 1 to functions compatible with Fuchsian groups
given by Theorem 2 requires a generalization of Lemma 5 to quasiconformal deformations of Fuchsian groups. This generalization (and even to appropriate
Kleinian groups with admissible in some sense collections of non-invariant components) is given in \cite{Kr1}.

\bigskip
{\small\em{ \leftline{Department of Mathematics, Bar-Ilan
University, 5290002 Ramat-Gan, Israel} \leftline{and
Department of Mathematics, University of Virginia,  Charlottesville, VA 22904-4137, USA}}

\end{document}